\documentclass[11pt,reqno]{amsart}
\usepackage{graphicx}
\usepackage{verbatim}
\usepackage{textcomp}
\usepackage{amssymb}
\usepackage{cite}
\usepackage{amsmath}
\usepackage{latexsym}
\usepackage{amscd}
\usepackage{amsthm}
\usepackage{mathrsfs}
\usepackage{xypic}
\usepackage{bm}
\usepackage{url}
\usepackage{hyperref}

\newtheorem{thm}{Theorem}[section]
\newtheorem{corr}[thm]{Corollary}
\newtheorem{lem}[thm]{Lemma}

\theoremstyle{definition}

\theoremstyle{remark}
\newtheorem{rem}{Remark}[section]
\numberwithin{equation}{section}
\begin{document}
\title[Rigidity of Riemannian manifolds]
{Rigidity of Riemannian manifolds with positive scalar curvature}

\author{Guangyue Huang}
\address{Department of Mathematics, Henan Normal
University, Xinxiang 453007, P.R. China} \email{hgy@henannu.edu.cn }

\thanks{The research of authors is supported by NSFC.}

\maketitle

\begin{abstract}
For the Bach-flat closed manifold with positive scalar curvature, we prove a rigidity result under a given inequality involving the Weyl curvature and the traceless Ricci curvature. Moveover, under an inequality involving $L^{\frac{n}{2}}$-norm of the Weyl curvature, the traceless Ricci curvature and the Yamabe invariant, we also provide a similar rigidity result. As an application, we obtain some rigidity results on $4$-dimensional manifolds.
\end{abstract}



\section{Introduction}

We recall that a Riemannian manifold $(M^{n},g)$ is a gradient shrinking Ricci soliton if there exists a smooth function $f$ such that
\begin{equation}\label{add1-Sec-1}
R_{ij}+f_{ij}=\lambda g_{ij}
\end{equation}
for some positive constant $\lambda$. In \cite{Catino2016}, relying on sharp algebraic curvature estimates, Catino proved some rigidity results for closed gradient shrinking Ricci soliton satisfying a $L^{\frac{n}{2}}$-pinching condition. From \eqref{add1-Sec-1}, it is easy to see that the Ricci curvature can be expressed by the Hessian of function $f$, in some sense. However, on a common Riemannian manifold $(M^{n},g)$, we can not find a function which is related to the Ricci curvature. It is nature to ask whether can one obtain some rigidity results under the $L^{\frac{n}{2}}$-pinching condition, which are analogous to those of Catino in \cite{Catino2016} on a given Riemannian manifold.

In order to study conformal relativity, R. Bach \cite{Bach21} in early
1920s' introduced the following Bach tensor:
\begin{equation}\label{1-Sec-1}
B_{ij}=\frac{1}{n-3}W_{ikjl,lk}+\frac{1}{n-2}W_{ikjl}R_{kl},
\end{equation}
where $n\geq4$, $W_{ijkl}$ denotes the Weyl curvature.
A metric $g$ is called Bach-flat if $B_{ij}=0$.
The aim of this paper is to achieve some rigidity results under the $L^{\frac{n}{2}}$-pinching condition on a given Riemannian manifold. In order to state our results, throughout this paper, we always denote $\mathring{R}_{ij}$ the traceless Ricci curvature.

Now, we can state our first result as follows:

\begin{thm}\label{thm1-1}
Let $(M^{4},g)$ be a closed Bach-flat manifold with
positive constant scalar curvature. If
\begin{equation}\label{1Th-1}\aligned
\Big|W+\frac{\sqrt{2}}{3}\mathring{{\rm Ric}} \mathbin{\bigcirc\mkern-15mu\wedge} g\Big|<\frac{2}{3\sqrt{3}}R,
\endaligned\end{equation}
then $M^{4}$ is isometric to a quotient of the round sphere $\mathbb{S}^{4}$.

\end{thm}

\begin{rem}\label{Rem-1}
Under the condition
\begin{equation}\label{1Rem-1}\aligned
|W|+|\mathring{R}_{ij}|<\frac{R}{4n},
\endaligned\end{equation}
Fang and Yuan, in \cite{Yuan2017}, have proved (see \cite[Theorem A]{Yuan2017}) the same result as Theorem \ref{thm1-1}.
It is easy to see from \eqref{1Th-1}
\begin{equation}\label{1Rem-2}\aligned
|W|^2+|\mathring{R}_{ij}|^2
\leq|W|^2+\frac{16}{9}|\mathring{R}_{ij}|^2<\frac{4 R^2}{27},
\endaligned\end{equation}
which gives
\begin{equation}\label{1Rem-3}\aligned
|W|+|\mathring{R}_{ij}|\leq\sqrt{2(|W|^2+|\mathring{R}_{ij}|^2)}<\frac{4 R}{3\sqrt{6}}.
\endaligned\end{equation}
Clearly, we have
\begin{equation}\label{1Rem-4}\aligned
\frac{4 R}{3\sqrt{6}}>\frac{R}{16}
\endaligned\end{equation}
and hence, for $n=4$, our Theorem \ref{thm1-1} generalizes Theorem A of Fang and Yuan \cite{Yuan2017}.

\end{rem}

The Yamabe invariant $Y(M,[g])$ associated to $(M^n,g)$ is defined by
\begin{equation}\label{2-Sec-28}\aligned
Y(M,[g])=&\inf\limits_{\tilde{g}\in[g]
}\frac{\int_M\tilde{R}\,dv_{\tilde{g}}}{(\int_M\,dv_{\tilde{g}})^{\frac{n-2}{n}}}\\
=&\frac{4(n-1)}{n-2}\inf\limits_{u\in W^{1,2}(M)}\frac{\int_M|\nabla
u|^2\,dv_g+\frac{n-2}{4(n-1)}\int_M
Ru^2\,dv_g}{(\int_M|u|^{\frac{2n}{n-2}}\,dv_g)^{\frac{n-2}{n}}},
\endaligned\end{equation}
where $[g]$ is the conformal class of the metric $g$. For closed manifolds, $Y(M,[g])$ is positive if and only if there exists a conformal metric in $[g]$ with everywhere positive scalar curvature. Therefore,
for any closed manifold with positive scalar curvature, from \eqref{2-Sec-28}, we obtain
\begin{equation}\label{2-Sec-29}\aligned
\frac{n-2}{4(n-1)}Y(M,[g])(\int_M|u|^{\frac{2n}{n-2}}\,dv_g)^{\frac{n-2}{n}}\leq&\int_M|\nabla
u|^2\,dv_g+\frac{n-2}{4(n-1)}\int_MRu^2\,dv_g.
\endaligned\end{equation}

By the aid of the Yamabe invariant given in \eqref{2-Sec-28}, we can prove the following result:

\begin{thm}\label{thm1-2}
Let $(M^{n},g)$ be a closed Bach-flat manifold with
positive constant scalar curvature. If either $4\leq n\leq5$ and
\begin{equation}\label{add-2Th-1}\aligned
\Big(\int_M\Big|W+\frac{\sqrt{n}}{\sqrt{8}(n-2)}\mathring{{\rm Ric}} \mathbin{\bigcirc\mkern-15mu\wedge} g\Big|^{\frac{n}{2}}\,dv_g\Big)^{\frac{2}{n}}<\frac{1}{4}\sqrt{\frac{n-2}{2(n-1)}}Y(M,[g]);
\endaligned\end{equation}
or $n=6$ and
\begin{equation}\label{add-2Th-2}\aligned
\Big(\int_M\Big|W+\frac{\sqrt{3}}{8}\mathring{{\rm Ric}} \mathbin{\bigcirc\mkern-15mu\wedge} g\Big|^{3}\,dv_g\Big)^{\frac{1}{3}}<\frac{1}{25}\sqrt{\frac{21}{10}}Y(M,[g]),
\endaligned\end{equation}
then $M^{n}$ is isometric to a quotient of the round sphere $\mathbb{S}^{n}$.

\end{thm}

For manifolds with harmonic curvature, we also obtain the following similar rigidity results:

\begin{thm}\label{thm1-3}
Let $M^{4}$ be a closed manifold with harmonic curvature.
If the scalar curvature is positive and
\begin{equation}\label{3Th-1}\aligned
\int_M\Big|W+\frac{1}{\sqrt{2}}\mathring{{\rm Ric}} \mathbin{\bigcirc\mkern-15mu\wedge} g\Big|^2\,dv_g<\frac{25}{486}Y^2(M,[g]),
\endaligned\end{equation}
then $M^{4}$ is isometric to a quotient of the round sphere $\mathbb{S}^{4}$.

\end{thm}

When $n=4$, we can get the following results:

\begin{corr}\label{2corr-1}
Let $(M^{4}, g)$ be a closed manifold. If it is Bach-flat and satisfies
\begin{equation}\label{2corrformula-1}
\int_M\Big(|W|^2+\frac{5}{4}|\mathring{R}_{ij}|^2\Big)\,dv_g
\leq\frac{1}{48}\int_MR^2\,dv_g,
\end{equation}
then it is isometric to a quotient of the round sphere $\mathbb{S}^{4}$.
If the curvature is harmonic and
\begin{equation}\label{2corrformula-2}
\int_M\Big(|W|^2+\frac{374}{81}|\mathring{R}_{ij}|^2\Big)\,dv_g
\leq\frac{25}{486}\int_MR^2\,dv_g,
\end{equation}
then it is also isometric to a quotient of the round sphere $\mathbb{S}^{4}$.

\end{corr}

\begin{rem}\label{2rem-1111}
The pinching conditions \eqref{2corrformula-1} and \eqref{2corrformula-2} are equivalent to
\begin{equation}\label{2rem-1}
\frac{13}{8}\int_M|W|^2\,dv_g+\frac{1}{12}\int_MR^2\,dv_g\leq20\pi^2\chi(M)
\end{equation}
and
\begin{equation}\label{2rem-2}
\frac{268}{81}\int_M|W|^2\,dv_g+\frac{1}{3}\int_MR^2\,dv_g
\leq\frac{5984}{81}\pi^2\chi(M),
\end{equation}
respectively. Here $\chi(M)$ is the Euler-Poincar\'{e} characteristic of $M^4$.

\end{rem}

\section{Some lemmas}

Recall that the Weyl curvature $W_{ijkl}$ is defined by
\begin{equation}\label{1-Sec-2}\aligned
W_{ijkl}=&R_{ijkl}-\frac{1}{n-2}(R_{ik}g_{jl}-R_{il}g_{jk}
+R_{jl}g_{ik}-R_{jk}g_{il})\\
&+\frac{R}{(n-1)(n-2)}(g_{ik}g_{jl}-g_{il}g_{jk})\\
=&R_{ijkl}-\frac{1}{n-2}(\mathring{R}_{ik}g_{jl}-\mathring{R}_{il}g_{jk}
+\mathring{R}_{jl}g_{ik}-\mathring{R}_{jk}g_{il})\\
&-\frac{R}{n(n-1)}(g_{ik}g_{jl}-g_{il}g_{jk}),
\endaligned
\end{equation}
where $R$ is the scalar curvature. Since the divergence of the Weyl curvature tensor is related to the Cotton tensor by
\begin{equation}\label{1-Sec-3}
-\frac{n-3}{n-2}C_{ijk}=W_{ijkl,l},
\end{equation}
where the Cotton tensor is given by
\begin{equation}\label{1-Sec-4}\aligned
C_{ijk}=&R_{kj,i}-R_{ki,j}-\frac{1}{2(n-1)}(R_{,i}g_{jk}-R_{,j}g_{ik})\\
=&\mathring{R}_{kj,i}-\mathring{R}_{ki,j}+\frac{n-2}{2n(n-1)}(R_{,i}g_{jk}-R_{,j}g_{ik}),
\endaligned
\end{equation}
the formula \eqref{1-Sec-1} can be written
as
\begin{equation}\label{1-Sec-5}
B_{ij}=\frac{1}{n-2}(C_{kij,k}+W_{ikjl}R^{kl}).
\end{equation}

As in \cite{Yuan2017}, for any $\theta\in \mathbb{R}$ and a symmetric 2-tensor $\varphi$, we introduce the following $\theta$-tensor:
\begin{equation}\label{2-Sec-1}
C_{\theta}^{\varphi}{}_{ijk}:=\varphi_{kj,i}-\theta\varphi_{ki,j}.
\end{equation}
We call a $\theta$-tensor $C_{\theta}^{\varphi}$ is the $\theta$-Codazzi tensor if $C_{\theta}^{\varphi}=0$. In particular, a Codazzi tensor is a special case of $\theta$-Codazzi tensor with $\theta=1$ and the
$\theta$-Codazzi tensor with $\theta=-1$ is called the anti-Codazzi tensor.

\begin{lem}\label{lem2-1}

Let $M^n$ be a closed manifold. Then for the $\theta$-tensor defined by \eqref{2-Sec-1}, we have
\begin{equation}\aligned\label{2-Sec-2}
\int_M[(\theta^2 +1)|\nabla \varphi|^2-|C_{\theta}^{\varphi}|^2]\,dv_g
=&2\theta\int_M\Big(|{\rm div}\varphi|^2+W_{ijkl}\mathring{\varphi}_{jl}\mathring{\varphi}_{ik}\\
&-\frac{n}{n-2}\mathring{R}_{ij}\mathring{\varphi}_{jk}
\mathring{\varphi}_{ki}-\frac{R}{n-1}|\mathring{\varphi}_{ij}|^2\Big)\,dv_g,
\endaligned\end{equation}
where $\mathring{\varphi}_{ij}=\varphi_{ij}-\frac{{\rm tr}\varphi}{n}g_{ij}$.
In particular, taking $\varphi_{ij}=R_{ij}$ in \eqref{2-Sec-2}, we obtain
\begin{equation}\aligned\label{add2-Sec-2}
\int_M|\nabla \mathring{R}_{ij}|^2\,dv_g\geq&\frac{2\theta}{\theta^2+1}\int_M\Big(
W_{ijkl}\mathring{R}_{jl}\mathring{R}_{ik}
-\frac{n}{n-2}\mathring{R}_{ij}\mathring{R}_{jk}
\mathring{R}_{ki}\\
&-\frac{1}{n-1}R|\mathring{R}_{ij}|^2+\frac{(n-2)^2}{4n^2}|\nabla R|^2\Big)\,dv_g.
\endaligned\end{equation}

\end{lem}

\proof For any $\theta\in \mathbb{R}$, we have
\begin{equation}\aligned\label{2-Sec-3}
\int_M|C_{\theta}^{\varphi}|^2\,dv_g=&\int_M|\varphi_{kj,i}-\theta\varphi_{ki,j}|^2\,dv_g\\
=&\int_M[(1+\theta^2)|\nabla \varphi|^2-2\theta\varphi_{kj,i}\varphi_{ki,j}]\,dv_g.
\endaligned\end{equation}
By the integration by parts, we have
\begin{equation}\aligned\label{2-Sec-4}
\int_M\varphi_{kj,i}\varphi_{ki,j}\,dv_g=&-\int_M\varphi_{kj,ij}\varphi_{ki}\,dv_g\\
=&-\int_M(\varphi_{kj,ji}+\varphi_{lj}R_{lkij}+\varphi_{kl}R_{ljij})\varphi_{ki}\,dv_g\\
=&\int_M(-\varphi_{kj,ji}+\varphi_{lj}R_{ijkl}-\varphi_{kl}R_{li})\varphi_{ki}\,dv_g\\
=&\int_M\Big[|{\rm div}\varphi|^2+R_{ijkl}\varphi_{jl}\varphi_{ik}
-\mathring{R}_{il}\mathring{\varphi}_{lk}\mathring{\varphi}_{ki}-\frac{2}{n}({\rm tr}\varphi)\mathring{R}_{ij}\mathring{\varphi}_{ij}\\
&-\frac{R}{n}\Big(|\mathring{\varphi}_{ij}|^2+\frac{({\rm tr}\varphi)^2}{n}\Big)\Big]\,dv_g.
\endaligned\end{equation}
Applying
\begin{equation}\aligned\label{2-Sec-5}
R_{ijkl}\varphi_{jl}\varphi_{ik}=&W_{ijkl}\varphi_{jl}\varphi_{ik}+\frac{2}{n-2}\Big[({\rm tr}\varphi)\mathring{R}_{ij}\varphi_{ij}-\mathring{R}_{ij}\varphi_{jk}\varphi_{ki}\Big]\\
&+\frac{R}{n(n-1)}[({\rm tr}\varphi)^2-|\varphi_{ij}|^2]\\
=&W_{ijkl}\mathring{\varphi}_{jl}\mathring{\varphi}_{ik}+\frac{2}{n-2}\Big[\frac{n-2}{n}({\rm tr}\varphi)\mathring{R}_{ij}\mathring{\varphi}_{ij}-\mathring{R}_{ij}\mathring{\varphi}_{jk}
\mathring{\varphi}_{ki}\Big]\\
&+\frac{R}{n(n-1)}\Big[\frac{n-1}{n}({\rm tr}\varphi)^2-|\mathring{\varphi}_{ij}|^2\Big]
\endaligned\end{equation}
in \eqref{2-Sec-4}, we obtain
\begin{equation}\aligned\label{2-Sec-6}
\int_M\varphi_{kj,i}\varphi_{ki,j}\,dv_g=&\int_M\Big[|{\rm div}\varphi|^2+W_{ijkl}\mathring{\varphi}_{jl}\mathring{\varphi}_{ik}+\frac{2}{n}({\rm tr}\varphi)\mathring{R}_{ij}\mathring{\varphi}_{ij}\\
&-\frac{n}{n-2}\mathring{R}_{ij}\mathring{\varphi}_{jk}
\mathring{\varphi}_{ki}-\frac{R}{n-1}|\mathring{\varphi}_{ij}|^2\Big]\,dv_g
\endaligned\end{equation}
and the desired \eqref{2-Sec-2} follows from by putting \eqref{2-Sec-6} into \eqref{2-Sec-3}.

Using the second Bianchi identity, we have $\mathring{R}_{ij,j}=\frac{n-2}{2n}R_{,i}$. Hence, \eqref{add2-Sec-2} follows by letting $\varphi_{ij}=R_{ij}$ in \eqref{2-Sec-2}. Thus, we conclude the proof of Lemma \ref{lem2-1}.
\endproof

\begin{lem}\label{lem2-2}

Let $M^n$ be a closed Bach-flat Riemannian manifold. Then
we have
\begin{equation}\aligned\label{2-Sec-8}
\int_M|\nabla \mathring{R}_{ij}|^2\,dv_g=&\int_M\Big(2W_{ijkl}\mathring{R}_{jl}\mathring{R}_{ik}
-\frac{n}{n-2}\mathring{R}_{ij}\mathring{R}_{jk}
\mathring{R}_{ki}\\
&-\frac{1}{n-1}R|\mathring{R}_{ij}|^2+\frac{(n-2)^2}{4n(n-1)}|\nabla R|^2\Big)\,dv_g.
\endaligned\end{equation}

\end{lem}

\proof Using the formula \eqref{1-Sec-2}, we can derive
\begin{equation}\label{2-Sec-9}\aligned
\mathring{R}_{kl}R_{ikjl}=&\mathring{R}_{kl}W_{ikjl}+\frac{1}{n-2}(|\mathring{R}_{ij}|^2g_{ij}
-2\mathring{R}_{ik}\mathring{R}_{jk})-\frac{1}{n(n-1)}R \mathring{R}_{ij},
\endaligned\end{equation}
which shows
\begin{equation}\label{2-Sec-10}\aligned
\mathring{R}_{kj,ik}=&\mathring{R}_{kj,ki}+\mathring{R}_{lj}R_{lkik}+\mathring{R}_{kl}R_{ljik}\\
=&\frac{n-2}{2n}R_{,ij}+\mathring{R}_{ik}\mathring{R}_{jk}+\frac{1}{n}R\mathring{R}_{ij}-\Big[\mathring{R}_{kl}W_{ikjl}
\\
&+\frac{1}{n-2}(|\mathring{R}_{ij}|^2g_{ij}-2\mathring{R}_{ik}\mathring{R}_{jk})
-\frac{1}{n(n-1)}R \mathring{R}_{ij}\Big]\\
=&\frac{n-2}{2n}R_{,ij}+\frac{n}{n-2}\mathring{R}_{ik}\mathring{R}_{jk}+\frac{1}{n-1}R\mathring{R}_{ij}
-\mathring{R}_{kl}W_{ikjl}-\frac{1}{n-2}|\mathring{R}_{ij}|^2g_{ij}.
\endaligned\end{equation}
Thus, from \eqref{1-Sec-4} and \eqref{2-Sec-10}, we have
\begin{equation}\label{2-Sec-11}\aligned
C_{kij,k}=&\Delta \mathring{R}_{ij}-\mathring{R}_{kj,ik}+\frac{n-2}{2n(n-1)}(g_{ij}\Delta R-R_{,ij})\\
=&\Delta \mathring{R}_{ij}-\Big(\frac{n-2}{2n}R_{,ij}+\frac{n}{n-2}\mathring{R}_{ik}\mathring{R}_{jk}
+\frac{1}{n-1}R\mathring{R}_{ij}-\mathring{R}_{kl}W_{ikjl}\\
&-\frac{1}{n-2}|\mathring{R}_{ij}|^2g_{ij}\Big)+\frac{n-2}{2n(n-1)}(g_{ij}\Delta R-R_{,ij})
\endaligned\end{equation}
and
\begin{equation}\label{2-Sec-12}\aligned
0=&(n-2)B_{ij}\mathring{R}_{ij}\\
=&C_{kij,k}\mathring{R}_{ij}+W_{ikjl}\mathring{R}_{ij}\mathring{R}_{kl}\\
=&\mathring{R}_{ij}\Delta \mathring{R}_{ij}-\frac{n}{n-2}\mathring{R}_{ij}\mathring{R}_{jk}\mathring{R}_{ki}
-\frac{1}{n-1}R|\mathring{R}_{ij}|^2\\
&+2W_{ikjl}\mathring{R}_{ij}\mathring{R}_{kl}-\frac{n-2}{2(n-1)}\mathring{R}_{ij}R_{,ij},
\endaligned\end{equation}
which gives
\begin{equation}\label{2-Sec-13}\aligned
\int_M|\nabla \mathring{R}_{ij}|^2\,dv_g=&-\int_M\mathring{R}_{ij}\Delta\mathring{R}_{ij}\,dv_g\\
=&\int_M\Big(2W_{ijkl}\mathring{R}_{jl}\mathring{R}_{ik}
-\frac{n}{n-2}\mathring{R}_{ij}\mathring{R}_{jk}
\mathring{R}_{ki}\\
&-\frac{1}{n-1}R|\mathring{R}_{ij}|^2+\frac{(n-2)^2}{4n(n-1)}|\nabla R|^2\Big)\,dv_g.
\endaligned\end{equation}
We complete the proof of Lemma \ref{lem2-2}.
\endproof

We introduce the following algebraic curvature estimate involving the Weyl curvature and traceless Ricci curvature:

\begin{lem}\label{add-Lemma22}
On every $n$-dimensional Riemannian manifold $(M^n,g)$, for any $\rho\in \mathbb{R}$, the following estimate holds
\begin{equation}\label{add2-Sec-21}\aligned
\Big|&-W_{ijkl}\mathring{R}_{jl}\mathring{R}_{ik}
+\frac{\rho}{n-2}\mathring{R}_{ij}\mathring{R}_{jk}\mathring{R}_{ki}\Big|\\
\leq&
\sqrt{\frac{n-2}{2(n-1)}}\Big(|W|^2+\frac{2\rho^2}{n(n-2)}|\mathring{R}_{ij}|^2
\Big)^{\frac{1}{2}}|\mathring{R}_{ij}|^2\\
=&\sqrt{\frac{n-2}{2(n-1)}}\Big|W+\frac{\rho}{\sqrt{2n}(n-2)}\mathring{{\rm Ric}} \mathbin{\bigcirc\mkern-15mu\wedge} g\Big||\mathring{R}_{ij}|^2.
\endaligned\end{equation}

\end{lem}

\proof We will prove Lemma \ref{add-Lemma22} only by some modifications for the proof of Proposition 2.1 of Catino in \cite{Catino2016}. Following  \cite[Lemma 4.7]{Bour2010},  we have
$$\aligned
(\mathring{{\rm Ric}} \mathbin{\bigcirc\mkern-15mu\wedge} g)_{ijkl}=&
\mathring{R}_{ik}g_{jl}-\mathring{R}_{il}g_{jk}+\mathring{R}_{jl}g_{ik}-\mathring{R}_{jk}g_{il},\\
(\mathring{{\rm Ric}} \mathbin{\bigcirc\mkern-15mu\wedge} \mathring{{\rm Ric}})_{ijkl}=&
2(\mathring{R}_{ik}\mathring{R}_{jl}-\mathring{R}_{il}\mathring{R}_{jk}).
\endaligned$$
It is easy to see
$$\aligned
W_{ijkl}\mathring{R}_{ik}\mathring{R}_{jl}=&\frac{1}{4}W_{ijkl}(\mathring{{\rm Ric}} \mathbin{\bigcirc\mkern-15mu\wedge} \mathring{{\rm Ric}})_{ijkl},\\
\mathring{R}_{ij}\mathring{R}_{jk}\mathring{R}_{ki}=&-\frac{1}{8}(\mathring{{\rm Ric}} \mathbin{\bigcirc\mkern-15mu\wedge} g)_{ijkl}(\mathring{{\rm Ric}} \mathbin{\bigcirc\mkern-15mu\wedge} \mathring{{\rm Ric}})_{ijkl},
\endaligned$$
which shows that
\begin{equation}\label{add2-Sec-23}
-W_{ijkl}\mathring{R}_{jl}\mathring{R}_{ik}
+\frac{\rho}{n-2}\mathring{R}_{ij}\mathring{R}_{jk}\mathring{R}_{ki}=-\frac{1}{4}
\Big(W+\frac{\rho}{2(n-2)}\mathring{{\rm Ric}} \mathbin{\bigcirc\mkern-15mu\wedge} g
\Big)_{ijkl}(\mathring{{\rm Ric}} \mathbin{\bigcirc\mkern-15mu\wedge} \mathring{{\rm Ric}})_{ijkl}.
\end{equation}
Let
\begin{equation}\label{add2-Sec-24}
T=\mathring{{\rm Ric}} \mathbin{\bigcirc\mkern-15mu\wedge} \mathring{{\rm Ric}}-U-V,
\end{equation}
where
$$\aligned
U_{ijkl}=&-\frac{2}{n(n-1)}|\mathring{R}_{ij}|^2(g \mathbin{\bigcirc\mkern-15mu\wedge} g)_{ijkl},\\
V_{ijkl}=&-\frac{2}{n-2}(\mathring{{\rm Ric}}^2 \mathbin{\bigcirc\mkern-15mu\wedge} g)_{ijkl}
+\frac{4}{n(n-2)}|\mathring{R}_{ij}|^2(g \mathbin{\bigcirc\mkern-15mu\wedge} g)_{ijkl}
\endaligned$$
with $(\mathring{{\rm Ric}}^2)_{ij}=\mathring{R}_{ik}\mathring{R}_{kj}$. Then $T$ is totally tracefree.
Taking the squard norm, we have
$$\aligned
|\mathring{{\rm Ric}} \mathbin{\bigcirc\mkern-15mu\wedge} \mathring{{\rm Ric}}|^2=&8|\mathring{R}_{ij}|^4-8|\mathring{{\rm Ric}}^2|^2,\\
|U|^2=&\frac{8}{n(n-1)}|\mathring{R}_{ij}|^4,\\
|V|^2=&\frac{16}{n-2}|\mathring{{\rm Ric}}^2|^2-\frac{16}{n(n-2)}|\mathring{R}_{ij}|^4
\endaligned$$
and
\begin{equation}\label{add2-Sec-25}
|T|^2+\frac{n}{2}|V|^2=|\mathring{{\rm Ric}} \mathbin{\bigcirc\mkern-15mu\wedge} \mathring{{\rm Ric}}|^2+\frac{n-2}{2}|V|^2-|U|^2=\frac{8(n-2)}{n-1}|\mathring{R}_{ij}|^4.
\end{equation}
Since both $W$ and $T$ are totally tracefree, using the Cauchy-Schwarz inequality, we obtain
\begin{equation}\label{add2-Sec-26}\aligned
\Big|\Big(W+&\frac{\rho}{2(n-2)}\mathring{{\rm Ric}} \mathbin{\bigcirc\mkern-15mu\wedge} g
\Big)_{ijkl}(\mathring{{\rm Ric}} \mathbin{\bigcirc\mkern-15mu\wedge} \mathring{{\rm Ric}})_{ijkl}\Big|^2\\
=&\Big|\Big(W+\frac{\rho}{2(n-2)}\mathring{{\rm Ric}} \mathbin{\bigcirc\mkern-15mu\wedge} g
\Big)_{ijkl}(T+V)_{ijkl}\Big|^2\\
=&\Big|\Big(W+\frac{\rho}{2(n-2)}\sqrt{\frac{2}{n}}\mathring{{\rm Ric}} \mathbin{\bigcirc\mkern-15mu\wedge} g
\Big)_{ijkl}(T+\sqrt{\frac{n}{2}}V)_{ijkl}\Big|^2\\
\leq&\Big|W+\frac{\rho}{2(n-2)}\sqrt{\frac{2}{n}}\mathring{{\rm Ric}} \mathbin{\bigcirc\mkern-15mu\wedge} g
\Big|^2\Big(|T|^2+\frac{n}{2}|V|^2\Big)\\
=&\frac{8(n-2)}{n-1}\Big(|W|^2+\frac{2\rho^2}{n(n-2)}|\mathring{R}_{ij}|^2\Big)|\mathring{R}_{ij}|^4.
\endaligned\end{equation}
Applying \eqref{add2-Sec-26} in \eqref{add2-Sec-23} and using
that
$$\Big|W+\frac{\rho}{\sqrt{2n}(n-2)}\mathring{{\rm Ric}} \mathbin{\bigcirc\mkern-15mu\wedge} g\Big|^2=|W|^2+\frac{2\rho^2}{n(n-2)}|\mathring{R}_{ij}|^2,$$
we obtain the desired estimate \eqref{add2-Sec-21}.

We complete the proof of Lemma \ref{add-Lemma22}.\endproof

Following the proof in \cite{Hebey1996}, for manifolds of which the metric $g$ is Einstein, we can prove the following (see \cite[Equation (5)]{Hebey1996})
\begin{equation}\label{2-Sec-19}\aligned
\Delta W_{ijkl}=\frac{2}{n}R W_{ijkl}-2(2W_{ipkq}W_{pjql}+\frac{1}{2}W_{klpq}W_{pqij}),
\endaligned\end{equation}
which gives
\begin{equation}\label{2-Sec-20}\aligned
\frac{1}{2}\Delta |W|^2=|\nabla W|^2+\frac{2}{n}R|W|^2-2(2W_{ijkl}W_{ipkq}W_{pjql}+\frac{1}{2}W_{ijkl}W_{klpq}W_{pqij}).
\endaligned\end{equation}

The following estimate comes from \cite[Lemma 2.5]{huang2017} (or see \cite{Catino2016,FuXiao2015}):

\begin{lem}\label{Lemma23}
On every $n$-dimensional Riemannian manifold $(M^n,g)$, there exists a positive constant $C(n)$
such that the following estimate holds
\begin{equation}\label{2-Sec-21}
2W_{ijkl}W_{ipkq}W_{pjql}+\frac{1}{2}W_{ijkl}W_{klpq}W_{pqij}\leq
C(n)|W|^3,
\end{equation}
where $C(n)$ is defined by
\begin{equation}\label{cn}
C(n)=
\begin{cases}\frac{\sqrt{6}}{4},\ \ \quad n=4\\
\frac{4\sqrt{10}}{15},  \quad n=5\\
\frac{\sqrt{70}}{2\sqrt{3}}, \ \quad n=6\\
\frac{5}{2}, \ \ \ \ \quad n\geq 7.
\end{cases}
\end{equation}

\end{lem}

Now we recall the following result which follows from the proof of \cite[Theorem 3.3]{Catino2016}:

\begin{lem}\label{2Lemma-4}

Let $C(n)$ be defined by \eqref{cn}. Then for the
Einstein manifold $(M^n,g)$, $n\geq 4$, with positive scalar
curvature, we have
\begin{equation}\label{2Compact9}\aligned
&\Bigg[\frac{n+1}{n-1}Y(M,[g])-\frac{8(n-1)}{n-2}C(n)
\Big(\int_M|W|^{\frac{n}{2}}\,dv_g\Big)^{\frac{2}{n}}\Bigg]
\int_M|\nabla|W||^2\,dv_g\\
&\qquad+\Bigg[\frac{2}{n}Y(M,[g])-2C(n)
\Big(\int_M|W|^{\frac{n}{2}}\,dv_g\Big)^{\frac{2}{n}}\Bigg]\int_MR|W|^2\,dv_g\leq0.
\endaligned\end{equation}

\end{lem}\label{3Lemma-1}

\proof Using \eqref{2-Sec-21}, we can deduce from \eqref{2-Sec-20}
\begin{equation}\label{2-Sec-22}\aligned
\frac{1}{2}\Delta |W|^2\geq&|\nabla W|^2+\frac{2}{n}R|W|^2-2C(n)|W|^3\\
\geq&2\Big(\frac{1}{n}R-C(n)|W|\Big)|W|^2.
\endaligned\end{equation}
Putting the following refined Kato inequality of any Einstein manifold
\begin{equation}\label{2-Sec-34}\aligned
|\nabla W|^2\geq\frac{n+1}{n-1}|\nabla |W||^2
\endaligned\end{equation}
into \eqref{2-Sec-22} yields
\begin{equation}\label{2-Sec-35}\aligned
\frac{1}{2}\Delta |W|^2\geq&\frac{n+1}{n-1}|\nabla |W||^2+\frac{2}{n}R|W|^2-2C(n)|W|^3.
\endaligned\end{equation}
Using the H\"{o}lder inequality and \eqref{2-Sec-29} with $u=|W|$, we get
\begin{equation}\label{2-Sec-36}\aligned
0\geq&\frac{n+1}{n-1}\int_M|\nabla |W||^2\,dv_g+\frac{2}{n}\int_MR|W|^2\,dv_g-2C(n)\int_M|W|^3\,dv_g\\
\geq&\frac{n+1}{n-1}\int_M|\nabla |W||^2\,dv_g+\frac{2}{n}\int_MR|W|^2\,dv_g\\
&-2C(n)\Big(\int_M|W|^{\frac{n}{2}}\,dv_g\Big)^{\frac{2}{n}}
\Big(\int_M|W|^{\frac{2n}{n-2}}\,dv_g\Big)^{\frac{n-2}{n}}\\
\geq&\frac{n+1}{n-1}\int_M|\nabla |W||^2\,dv_g+\frac{2}{n}\int_MR|W|^2\,dv_g\\
&-2\frac{C(n)}{Y(M,[g])}\Big(\int_M|W|^{\frac{n}{2}}\,dv_g\Big)^{\frac{2}{n}}
\Big(\frac{4(n-1)}{n-2}\int_M|\nabla|W||^2\,dv_g+\int_MR|\nabla W |^2\,dv_g\Big)\\
=&\Bigg[\frac{n+1}{n-1}Y(M,[g])-\frac{8(n-1)}{n-2}C(n)
\Big(\int_M|W|^{\frac{n}{2}}\,dv_g\Big)^{\frac{2}{n}}\Bigg]
\int_M|\nabla|W||^2\,dv_g\\
&\qquad+\Bigg[\frac{2}{n}Y(M,[g])-2C(n)
\Big(\int_M|W|^{\frac{n}{2}}\,dv_g\Big)^{\frac{2}{n}}\Bigg]\int_MR|W|^2\,dv_g,
\endaligned\end{equation}
and the proof of Lemma \ref{2Lemma-4} is completed. \endproof

\section{Proof of Theorems}

\subsection{Proof of Theorem \ref{thm1-1}}

Now, with the help of Lemma \ref{lem2-1} and Lemma \ref{lem2-2}, we will complete the proof of Theorem \ref{thm1-1}.

Combining \eqref{add2-Sec-2} with \eqref{2-Sec-8}, we obtain
\begin{equation}\label{2-Sec-14}\aligned
0\geq&\int_M\Big[2\Big(-W_{ijkl}\mathring{R}_{jl}\mathring{R}_{ik}
+\frac{n(\theta-1)^2}{2(\theta^2-\theta+1)}\frac{1}{n-2}\mathring{R}_{ij}\mathring{R}_{jk}
\mathring{R}_{ki}\Big)\\
&+\frac{(\theta-1)^2}{(n-1)(\theta^2-\theta+1)}R|\mathring{R}_{ij}|^2\Big]\,dv_g.
\endaligned\end{equation}
Since for any $\theta$, we have $\theta^2-\theta+1>0$,
using \eqref{add2-Sec-21}, we can deduce from \eqref{2-Sec-14}
\begin{equation}\label{2-Sec-17}\aligned
0\geq&\int_M\Big[-\sqrt{\frac{2(n-2)}{n-1}}\Big|W+\frac{n(\theta-1)^2}{2\sqrt{2n}(n-2)(\theta^2-\theta+1)}\mathring{{\rm Ric}} \mathbin{\bigcirc\mkern-15mu\wedge} g\Big|\\
&+\frac{(\theta-1)^2}{(n-1)(\theta^2-\theta+1)}R\Big]|\mathring{R}_{ij}|^2\,dv_g\\
\geq&0
\endaligned\end{equation}
under the condition
\begin{equation}\label{2-Sec-18}\aligned
\sqrt{\frac{2(n-2)}{n-1}}\Big|\frac{\theta^2-\theta+1}{(\theta-1)^2}W
+\frac{n}{2\sqrt{2n}(n-2)}\mathring{{\rm Ric}} \mathbin{\bigcirc\mkern-15mu\wedge} g\Big|<\frac{1}{n-1}R
\endaligned\end{equation}
with a given $\theta$ ($\theta\neq1$). Minimizing the  coefficient  of  $W$ with respect to the function $\theta$ by taking
\begin{equation}\label{2-Sec-25}
\theta=-1,
\end{equation}
we obtain that if
\begin{equation}\label{add-1Th-1}\aligned
\Big|W+\frac{\sqrt{2n }}{3(n-2)}\mathring{{\rm Ric}} \mathbin{\bigcirc\mkern-15mu\wedge} g\Big|<\frac{4}{3\sqrt{2(n-1)(n-2)}}R
\endaligned\end{equation}
holds, then \eqref{2-Sec-17} shows
$\mathring{R}_{ij}=0$ and the metric $g$ must be Einstein.

On the other hand, integrating \eqref{2-Sec-22}, we have
\begin{equation}\label{2-Sec-23}\aligned
\int_M\Big(\frac{1}{n}R-C(n)|W|\Big)|W|^2\,dv_g\leq0.
\endaligned\end{equation}
Under the assumption \eqref{1Th-1}, we have that the metric $g$ is Einstein and the condition \eqref{add-1Th-1} becomes
\begin{equation}\label{2-Sec-26}\aligned
|W|<\frac{4}{3\sqrt{2(n-1)(n-2)}}R.
\endaligned\end{equation}
In this case, we have
\begin{equation}\label{2-Sec-26}\aligned
C(n)|W|<\frac{4C(n)}{3\sqrt{2(n-1)(n-2)}}R<\frac{1}{n}R,
\endaligned\end{equation}
for $n=4$, and \eqref{2-Sec-23} yields
\begin{equation}\label{2-Sec-27}\aligned
0\leq\int_M\Big(\frac{1}{n}R-C(n)|W|\Big)|W|^2\,dv_g\leq0.
\endaligned\end{equation}
This shows that $W_{ijkl}=0$ and $M^{4}$ is isometric to a quotient of the round sphere $\mathbb{S}^{4}$.

We complete the proof of Theorem \ref{thm1-1}.

\subsection{Proof of Theorem \ref{thm1-2}}

Replacing $u$ in \eqref{2-Sec-29} with $|\mathring{R}_{ij}|$ yields
\begin{equation}\label{2-Sec-30}\aligned
\frac{n-2}{4(n-1)}Y(M,[g])&(\int_M|\mathring{R}_{ij}|^{\frac{2n}{n-2}}\,dv_g)^{\frac{n-2}{n}}\\
\leq&\int_M|\nabla
|\mathring{R}_{ij}||^2\,dv_g+\frac{n-2}{4(n-1)}\int_MR|\mathring{R}_{ij}|^2\,dv_g.
\endaligned\end{equation}
On the other hand, using the Kato inequality, \eqref{2-Sec-8} and \eqref{add2-Sec-21}, we obtain
\begin{equation}\aligned\label{2-Sec-31}
\int_M|\nabla
|\mathring{R}_{ij}||^2\,dv_g\leq&\int_M|\nabla \mathring{R}_{ij}|^2\,dv_g\\
=&\int_M\Big[2W_{ijkl}\mathring{R}_{jl}\mathring{R}_{ik}
-\frac{n}{n-2}\mathring{R}_{ij}\mathring{R}_{jk}
\mathring{R}_{ki}-\frac{1}{n-1}R|\mathring{R}_{ij}|^2\Big]\,dv_g\\
\leq&\int_M\Big[\sqrt{\frac{2(n-2)}{n-1}}\Big|W+\frac{\sqrt{n}}{\sqrt{8}(n-2)}\mathring{{\rm Ric}} \mathbin{\bigcirc\mkern-15mu\wedge} g\Big|-\frac{1}{n-1}R\Big]|\mathring{R}_{ij}|^2\,dv_g.
\endaligned\end{equation}
Inserting \eqref{2-Sec-31} into \eqref{2-Sec-30} gives
\begin{equation}\label{2-Sec-32}\aligned
\frac{n-2}{4(n-1)}Y(M,[g])&\Big(\int_M|\mathring{R}_{ij}|^{\frac{2n}{n-2}}\,dv_g\Big)^{\frac{n-2}{n}}\\
\leq&\int_M\Big[\sqrt{\frac{2(n-2)}{n-1}}\Big|W+\frac{\sqrt{n}}{\sqrt{8}(n-2)}\mathring{{\rm Ric}} \mathbin{\bigcirc\mkern-15mu\wedge} g\Big|\\
&+\frac{n-6}{4(n-1)}R\Big]|\mathring{R}_{ij}|^2\,dv_g.
\endaligned\end{equation}
Hence, for $4\leq n\leq6$, using H\"{o}lder inequality for \eqref{2-Sec-32} yields
\begin{equation}\label{2-Sec-37}\aligned
\Bigg[\frac{n-2}{4(n-1)}Y(M,[g])&-\sqrt{\frac{2(n-2)}{n-1}}\Big(\int_M\Big|W+\frac{\sqrt{n}}{\sqrt{8}(n-2)}\mathring{{\rm Ric}} \mathbin{\bigcirc\mkern-15mu\wedge} g\Big|^{\frac{n}{2}}\,dv_g\Big)^{\frac{2}{n}}\Bigg]\\
&\times\Big(\int_M|\mathring{R}_{ij}|^{\frac{2n}{n-2}}\,dv_g\Big)^{\frac{n-2}{n}}
\leq0,
\endaligned\end{equation}
which, under the condition
\begin{equation}\label{2-Sec-38}\aligned
\Big(\int_M\Big|W+\frac{\sqrt{n}}{\sqrt{8}(n-2)}\mathring{{\rm Ric}} \mathbin{\bigcirc\mkern-15mu\wedge} g\Big|^{\frac{n}{2}}\,dv_g\Big)^{\frac{2}{n}}
<\frac{1}{4}\sqrt{\frac{n-2}{2(n-1)}}Y(M,[g])
\endaligned\end{equation}
\eqref{2-Sec-37} shows that the metric $g$ must be Einstein. Since the metric $g$ is Einstein, then \eqref{2-Sec-38} becomes
\begin{equation}\label{2-Sec-40}\aligned
\Big(\int_M|W|^{\frac{n}{2}}\,dv_g\Big)^{\frac{2}{n}}
<\frac{1}{4}\sqrt{\frac{n-2}{2(n-1)}}Y(M,[g]).
\endaligned\end{equation}
Noticing for $4\leq n\leq5$,
\begin{equation}\label{add2-Sec-41}\aligned
\frac{1}{4}\sqrt{\frac{n-2}{2(n-1)}}=\min\Bigg\{\frac{1}{4}\sqrt{\frac{n-2}{2(n-1)}},\ \frac{(n+1)(n-2)}{8(n-1)^2C(n)},\ \frac{1}{n\,C(n)}\Bigg\},
\endaligned\end{equation}
and  for $n=6$,
\begin{equation}\label{add2-Sec-41}\aligned
\frac{(n+1)(n-2)}{8(n-1)^2C(n)}=\min\Bigg\{\frac{1}{4}\sqrt{\frac{n-2}{2(n-1)}},\ \frac{(n+1)(n-2)}{8(n-1)^2C(n)},\ \frac{1}{n\,C(n)}\Bigg\},
\endaligned\end{equation}
hence under the assumption \eqref{add-2Th-1} or \eqref{add-2Th-2}, we have
\begin{equation}\label{2-Sec-43}\aligned
0\leq&\Bigg[\frac{n+1}{n-1}Y(M,[g])-\frac{8(n-1)}{n-2}C(n)
\Big(\int_M|W|^{\frac{n}{2}}\,dv_g\Big)^{\frac{2}{n}}\Bigg]
\int_M|\nabla|W||^2\,dv_g\\
&+\Bigg[\frac{2}{n}Y(M,[g])-2C(n)
\Big(\int_M|W|^{\frac{n}{2}}\,dv_g\Big)^{\frac{2}{n}}\Bigg]\int_MR|W|^2\,dv_g\leq0,
\endaligned\end{equation}
which shows that $W_{ijkl}=0$ and $M^{n}$ is isometric to a quotient of the round sphere $\mathbb{S}^{n}$.

\subsection{Proof of Theorem \ref{thm1-3}}

If the curvature of $(M^n,g)$ is harmonic, then we have
\begin{equation}\label{ma-Rem-1}
R_{ki,j}-R_{kj,i}=R_{ijkl,l}=0,
\end{equation}
This shows that the Ricci curvature is Codazzi and
\begin{equation}\label{1-Rem-2}
R_{,i}=R_{kk,i}=R_{ik,k}=\frac{1}{2}R_{,i},\end{equation}
which shows the scalar curvature is constant. Therefore, \eqref{add2-Sec-2} becomes
\begin{equation}\aligned\label{3-Sec-1}
\int_M|\nabla \mathring{R}_{ij}|^2\,dv_g=&\int_M\Big[
W_{ijkl}\mathring{R}_{jl}\mathring{R}_{ik}
-\frac{n}{n-2}\mathring{R}_{ij}\mathring{R}_{jk}
\mathring{R}_{ki}\\
&-\frac{1}{n-1}R|\mathring{R}_{ij}|^2\Big]\,dv_g.
\endaligned\end{equation}

As first observed by Bourguignon \cite{Bourguignon1990}, the
Codazzi tensor $\mathring{R}_{ij}$ satisfies the following sharp inequality
(for a proof, see for instance \cite{Hebey1996}):
\begin{equation}\label{3-Sec-5}
|\nabla \mathring{R}_{ij}|^2\geq \frac{n+2}{n}|\nabla|\mathring{R}_{ij}||^2.
\end{equation}
It follows from \eqref{3-Sec-1} that
\begin{equation}\label{3-Sec-7}\aligned
0\geq&\frac{n+2}{n}\int_M|\nabla |\mathring{R}_{ij}||^2\,dv_g
+\int_M\Big[
-W_{ijkl}\mathring{R}_{jl}\mathring{R}_{ik}
+\frac{n}{n-2}\mathring{R}_{ij}\mathring{R}_{jk}
\mathring{R}_{ki}\\
&+\frac{1}{n-1}R|\mathring{R}_{ij}|^2\Big]\,dv_g,
\endaligned\end{equation}
which combining with \eqref{2-Sec-30} gives
\begin{equation}\label{3-Sec-8}\aligned
0\geq&\frac{n+2}{n}\Bigg[\frac{n-2}{4(n-1)}Y(M,[g])
\Big(\int_M|\mathring{R}_{ij}|^{\frac{2n}{n-2}}\,dv_g\Big)^{\frac{n-2}{n}}
-\frac{n-2}{4(n-1)}\int_MR|\mathring{R}_{ij}|^2\,dv_g\Bigg]\\
&-\int_M\sqrt{\frac{n-2}{2(n-1)}}\Big|W+\frac{\sqrt{n}}{\sqrt{2}(n-2)}\mathring{{\rm Ric}} \mathbin{\bigcirc\mkern-15mu\wedge} g\Big||\mathring{R}_{ij}|^2\,dv_g+\frac{1}{n-1}\int_MR|\mathring{R}_{ij}|^2\,dv_g\\
\geq&\Bigg[\frac{n^2-4}{4n(n-1)}Y(M,[g])
-\sqrt{\frac{n-2}{2(n-1)}}\Big(\int_M\Big|W+\frac{\sqrt{n}}{\sqrt{2}(n-2)}\mathring{{\rm Ric}} \mathbin{\bigcirc\mkern-15mu\wedge} g\Big|^{\frac{n}{2}}\Big)^{\frac{2}{n}}\Bigg]\\
&\times\Big(\int_M|\mathring{R}_{ij}|^{\frac{2n}{n-2}}\,dv_g\Big)^{\frac{n-2}{n}}
-\frac{n^2-4n-4}{4n(n-1)}\int_MR|\mathring{R}_{ij}|^2\,dv_g\\
\geq&0.
\endaligned\end{equation}
under \eqref{3Th-1} for $n=4$.
Hence, we obtain that the metric $g$ must be Einstein and hence \eqref{3Th-1} becomes
\begin{equation}\label{3-Sec-9}\aligned
\Big(\int_M|W|^{\frac{n}{2}}\,dv_g\Big)^{\frac{2}{n}}
<\frac{n+2}{2n}\sqrt{\frac{n-2}{2(n-1)}}Y(M,[g]).
\endaligned\end{equation}
Noticing for $n=4$,
\begin{equation}\label{3-Sec-10}\aligned
\frac{(n+1)(n-2)}{8(n-1)^2C(n)}=\min\Bigg\{\frac{n+2}{2n}\sqrt{\frac{n-2}{2(n-1)}},\ \frac{(n+1)(n-2)}{8(n-1)^2C(n)},\ \frac{1}{n\,C(n)}\Bigg\},
\endaligned\end{equation}
hence, under the condition \eqref{3Th-1}, we have that \eqref{3-Sec-8} holds and
\begin{equation}\label{3-Sec-12}\aligned
0\leq&\Bigg[\frac{n+1}{n-1}Y(M,[g])-\frac{8(n-1)}{n-2}C(n)
\Big(\int_M|W|^{\frac{n}{2}}\,dv_g\Big)^{\frac{2}{n}}\Bigg]
\int_M|\nabla|W||^2\,dv_g\\
&+\Bigg[\frac{2}{n}Y(M,[g])-2C(n)
\Big(\int_M|W|^{\frac{n}{2}}\,dv_g\Big)^{\frac{2}{n}}\Bigg]\int_MR|W|^2\,dv_g\leq0
\endaligned\end{equation}
which shows $W_{ijkl}=0$ and $M^{4}$ is isometric to a quotient of the round sphere $\mathbb{S}^{4}$.

We complete the proof of Theorem \ref{thm1-3}.

\section{Some proof for four dimensional manifolds}

When $n=4$, it has been prove by Catino (see \cite[Lemma 4.1]{Catino2016}) that

\begin{lem}\label{3-Lemma-1}
Let $(M^4,g)$ be a closed  manifold. Then
\begin{equation}\label{4-Sec-1}\aligned
Y^2(M,[g])\geq\int_M(R^2-12|\mathring{R}_{ij}|^2)\,dv_g,
\endaligned\end{equation}
with the inequlity is strict unless $(M^4,g)$  is conformally Einstein.

\end{lem}

Since \eqref{add-2Th-1} can be written as
\begin{equation}\label{4-Sec-2}\aligned
\int_M(|W|^2+|\mathring{R}_{ij}|^2)\,dv_g<\frac{1}{48}Y^2(M,[g]).
\endaligned\end{equation}
Using \eqref{4-Sec-1}, it is easy to see
\begin{equation}\label{4-Sec-3}\aligned
\int_M(|W^2+|\mathring{R}_{ij}|^2)\,dv_g-\frac{1}{48}Y^2(M,[g])
\leq\int_M\Big(|W^2+\frac{5}{4}|\mathring{R}_{ij}|^2-\frac{1}{48}R^2\Big)\,dv_g.
\endaligned\end{equation}
Moreover, \eqref{3Th-1} can be written as
\begin{equation}\label{4-Sec-4}\aligned
\int_M(|W^2+4|\mathring{R}_{ij}|^2)\,dv_g<\frac{25}{486}Y^2(M,[g]),
\endaligned\end{equation}
which combining with \eqref{3Th-1} gives
\begin{equation}\label{4-Sec-5}\aligned
\int_M(|W^2+4|\mathring{R}_{ij}|^2)\,dv_g-\frac{25}{486}Y^2(M,[g])\leq
\int_M\Big(|W^2+\frac{374}{81}|\mathring{R}_{ij}|^2-\frac{25}{486}R^2\Big)\,dv_g.
\endaligned\end{equation}
Hence, Corollary \ref{2corr-1} follows from \eqref{4-Sec-3} immediately.

When $n=4$, the following Chern-Gauss-Bonnet formula (see \cite[Equation 6.31]{Besse2008})
\begin{equation}\label{4-Sec-6}\aligned
\int_M\Big(|W^2-2|\mathring{R}_{ij}|^2+\frac{1}{6}R^2\Big)\,dv_g=32\pi^2\chi(M)
\endaligned\end{equation}
is well-known. This shows
\begin{equation}\label{4-Sec-7}\aligned
\int_M|\mathring{R}_{ij}|^2\,dv_g=\int_M\Big(\frac{1}{2}|W|^2+\frac{1}{12}R^2\Big)\,dv_g-16\pi^2\chi(M).
\endaligned\end{equation}
Therefore, \eqref{2rem-1} and \eqref{2rem-2} follow from inserting \eqref{4-Sec-7} into \eqref{2corrformula-1}
and \eqref{2corrformula-2}, respectively. This proves Remark \ref{2rem-1111}.


\bibliographystyle{Plain}

\end{document}